\documentclass[psamsfonts]{amsart}

\usepackage{amssymb,amsfonts}
\usepackage[all,arc]{xy}
\usepackage{enumerate}
\usepackage{mathrsfs}
\usepackage{graphicx}
\usepackage{tikz}
\usepackage{amscd}
\usepackage{amssymb}
\usepackage{xy}
\usepackage{lscape}
\usepackage{setspace}
\usepackage{amsthm}
\usepackage{mathtools}
\usetikzlibrary{matrix,arrows}

\newtheorem{thm}{Theorem}[section]

\newtheorem{lem}[thm]{Lemma}

\theoremstyle{definition}

\theoremstyle{remark}

\newcommand{\NN}{\mathbb N}
\newcommand{\ZZ}{\mathbb Z}

\makeatletter
\let\c@equation\c@thm
\makeatother
\numberwithin{equation}{section}

\title{Integral String Lie Algebra Structure of Spheres}

\author{Felicia Tabing}

\date{}

\begin{document}

\begin{abstract}

\indent Chas and Sullivan introduced string homology in \cite{CS}, which is the equivariant homology of the loop space with the $S^1$ action on loops by rotation. Craig Westerland computed the string homology for spheres with coefficients in $\mathbb{Z} /2\mathbb{Z}$ \cite{We} and Somnath Basu computed the string homology and string bracket for spheres over rational coefficients and found that the bracket is trivial in his dissertation \cite{Ba}. In this paper, we compute string homology and the string bracket for spheres with integer coefficients, treating the odd- and even-dimensional cases separately. We use the Gysin sequence and Leray-Serre spectral sequence for our computations. We find that over the integers, the string Lie algebra bracket structure is not always zero as Basu found.  The string bracket turns out to be non-zero on torsion elements coming from string homology. \\

\end{abstract}

\maketitle

\tableofcontents

\section{Introduction}

\indent The term \emph{String Topology} came from the paper of the same name by Moira Chas and Dennis Sullivan in 1999, where the authors introduced various algebraic structures that arose from the homology of the free loop space that came out generalizing the Lie algebra structure that William M. Goldman described by the intersection and concatenation of loops on surfaces \cite{Go}.\\
\indent In this paper, we compute the integral string homology and string bracket structure of spheres, where some torsion phenomena appear. In our computations, we use the Leray-Serre spectral sequence, and the Gysin exact sequence.  Craig Westerland computed this over $\mathbb{Z}_2$ using a spectral sequence and Somnath Basu computed this over a field  $\mathbb{Q}$ via minimal models.

\section{String Homology Lie Algebra Preliminaries}

\indent We describe the basic algebraic structures appearing in the homology and equivariant homology of the free loop space of an $n$ dimensional manifold $M$, $LM=Map(S^1, M)$ as described by Chas and Sullivan \cite{CS}.  We will consider homology and cohomology with integer coefficients, unless otherwise stated. We denote the usual homology of the free loop space of $M$ as $H_*(LM)$ and equivariant homology will be denoted by $H_*^{S^1}(LM)$.  For the rest of this paper we will drop the $LM$ from the homology notation if it is clear from the context which space $M$ we are referring to. \\
\indent Chas and Sullivan defined the \emph{loop product} on $\mathbb{H}_*(LM):=H_{*+d}(LM)$,
\begin{equation*}
-\bullet -= \mathbb{H}_*(LM) \otimes\mathbb{H}_*(LM) \rightarrow \mathbb{H}_{*}(LM),
\end{equation*}
which was defined on the chain level by intersection and concatenation of families of loops.
\emph{String homology} is the equivariant homology of the free loop space with respect to the action of rotation of loops, $H_*^{S^1}(LM)$.  The fibration,
\begin{align*}
S^2 \rightarrow LM \times ES^1 \rightarrow LM \times_{S^1} ES^1,
\end{align*}
induces a long exact sequence on homology, the Gysin sequence from which Chas and Sullivan used to describe a Lie bracket on $H_*^{S^1}(LM)$,
\begin{align*}
\cdots \rightarrow \mathbb{H}_{*-d}(LM) \xrightarrow{e} H_*^{S^1}(LM) \xrightarrow{\cap} H_{*-2}^{S^1}(LM) \xrightarrow{M} \mathbb{H}_{*-d-1}(LM) \rightarrow \cdots,
\end{align*}
where $e$ and $M$ are informally called the "erasing map" and "marking map," respectively. There is a \emph{Batalin-Vilkovisky operator} denoted by $\Delta$, which comes from the natural action given by rotation of loops,
\begin{align*}
\rho:S^1\times LM \rightarrow LM
\end{align*} 
given by $\rho(t,\gamma)(s)=\gamma(s+t)$. This action defines a degree one operator on loop homology: $\Delta:\mathbb{H}_*(LM)\rightarrow \mathbb{H}_{*+1}(LM)$ given by $\delta(\alpha)=\rho_*([S^1]\otimes \alpha)$ for $\alpha \in H_k(LM)$. Note that $e\circ M=0$ by exactness, and $\Delta = M \circ e$. \\
\indent For two classes $\alpha , \beta \in H_*^{S^1}(LM)$, the \emph{string bracket} is defined by
\begin{align*}
[\alpha , \beta ] = (-1)^{| \alpha | -d }e(M(\alpha) \bullet M(\beta) ), 
\end{align*}
where $\bullet$ is the loop product.
\thm (Chas-Sullivan) $(H_*^{S^1}(LM),[-,-])$ is a graded Lie algebra, with Lie bracket of degree $2-d$ \cite{CS}.\\
\normalfont
More precisely, our bracket is a map:
\begin{align*}
[-,-]: H_i^{S^1}(LM)\times H_j^{S^1}(LM) \rightarrow H_{i+j+2-d}^{S^1}(LM).
\end{align*}
\indent In the following sections, we compute the $H_*^{S^1}(LS^n)$ for all $n\in \NN$ and we compute the structure of the string bracket. 
\pagebreak

\section{String Homology and String Bracket Structure}

\thm \label{thm:stringhomology1} The string homology structure of spheres, $H^{S^1}_*(LS^n)$:
\begin{enumerate}
\item  for $n=1$,
\begin{eqnarray}
& H^{S^1}_0(LS^1) &\cong \bigoplus_{n\in \ZZ} \ZZ \nonumber \\
& H^{S^1}_1(LS^1) \cong  H^{S^1}_{2i+1}(LS^1) &\cong \bigoplus_{n\in \ZZ -\{0\}} \ZZ_n \oplus \ZZ , \hspace*{.5cm} i\geq 0\nonumber \\
& H^{S^1}_2(LS^1) \cong H^{S^1}_{2i}(LS^1) &\cong \ZZ , \hspace*{.5cm} i\geq 1\nonumber 
\end{eqnarray}\\
\item for $n=2k+1$,
\begin{equation} 
H^{S^1}_{2i}(LS^n) \cong \left\{
\begin{array}{rl}
\ZZ  & \text{if } (n-1)\nmid 2i,\\
\ZZ \oplus \ZZ & \text{if } (n-1)\vert 2i. \nonumber
\end{array} \right.
\end{equation}
for $1\leq i \leq \frac{(k+1)(n-1)-2}{2}$ and 
\begin{equation}
H_{2i+1}^{S^1}(LS^n)\cong t_k  \nonumber
\end{equation}
for $\frac{k(n-1)}{2}\leq i \leq \frac{(k+1)(n-1)-2}{2}$, where $t_k$ is a torsion group of order $k!$. All other $j$ that does not fall into the above categories, we have that $H_i^{S^1}(LS^n) \cong H_{i-2}^{S^1}(LS^n)$
\item for $n=2k$,
\begin{align*}
H_{2j+1}^{S^1}(LS^n)\cong H_{k(2n-2)+1}^{S^1}(LS^n)& \cong 0\\
H_{2j}^{S^1}(LS^n)\cong H_{k(2n-2)+n-1}^{S^1}(LS^n) & \cong \ZZ\\
H_{k(2n-2)+n}^{S^1}(LS^n) & \cong \ZZ_2^{k} \oplus C_k \oplus \ZZ \oplus \ZZ .
\end{align*}
for $j\leq \frac{n}{2}-1$ and all $k\in \ZZ\cup \{0\}$.  We have $H_i^{S^1}(LS^n)\cong H_{i-2}^{S^1}(LS^n)$ for $k(2n-2)+1 \leq i \leq k(2n-2)+n-1$ and $C_k$ is a torsion group of order  $\prod_{i=1}^{k-1}(2i+1)$.
\end{enumerate}
\pagebreak
\normalfont
\thm  \label{thm:stringbracket}The String Bracket structure of spheres, 
\begin{align*}
[-,-]: H_i^{S^1}(LS^n) \otimes H_j^{S^1}(LS^n) \rightarrow H_{i+j+2-n}^{S^1}(LS^n)
\end{align*}
is given as follows:
\begin{enumerate}
\item for $n=1$, the string bracket is only nontrivial on generators of degree zero.  For ${e(a \otimes x^n)}$, ${e(a \otimes x^m})$ in ${H^{S^1}_0(LS^1)}$,
\begin{align*}
[e(a \otimes x^n),e( a \otimes x^m)]=-nm(1 \otimes x^{n+m}),
\end{align*}
So ${[a \otimes x^n, a \otimes x^m]=0}$ if $n+m \neq 0$ and $n+m$ divides $nm$. If $n+m=0$ then  ${[a \otimes x^n, a \otimes x^m]=nm(1\otimes 1)}$.
The bracket is only nontrivial for the torsion elements. 
\item \label{thm:stringbracket2} for $n=2k+1$,  the only possible non-zero bracket is of the generators $e(a\otimes u^i)\in H_{2i}^{S^1}(LS^n)$ for $2i$ divisible by $n-1$.
\begin{align*}
[e(a\otimes u^i), e(a\otimes u^j)]=ije(1\otimes u^{i+j-2}),
\end{align*}
where $e(1\otimes u^{i+j-2})$ is a generator of $\ZZ_{i+j-1}$, so the bracket is only zero when $i+j-1$ divides $ij$.

\item for $n=2k$, the string bracket is always zero except on the generators $e(bv^j)$ of $H_{j(2n-2)+n-1}^{S^1}(LS^n)$.
\begin{align*}
[e(bv^k),e(bv^l)]=-(4kl+2k+2l+1)e(v^{k+l}),
\end{align*}
and since $e(v^{k+l})$ has order $2(k+l)+1$ it is not always zero.
\end{enumerate}
\normalfont
\section{The Gysin Sequence for the Circle Bundle}
We consider the even and odd spheres as separate cases in computing string homology.  In both cases, we use the Gysin sequence of the circle bundle,
\begin{align} \label{eqn:oddcirclebundle} 
 S^1 \rightarrow LS^n \times ES^1 \rightarrow LS^n \times_{S^1} ES^1. 
\end{align}
\subsection{String Homology for Odd Spheres}
\ \\
\indent For $n$ odd, the BV-operator acts on the generators of  $$\mathbb{H}_*(LS^n)=\Lambda[a]\otimes \ZZ [u]$$ as follows, \cite{Me}:
\begin{align*}
\Delta(a\otimes u^i)= & i(1 \otimes u^{i-1})\\
\Delta (1 \otimes u^i)= & 0.
\end{align*}
\normalfont
Consider the bottom of the Gysin sequence. Let $H^{S^1}_i$ denote $H^{S^1}_i(LS^n)$ and $\mathbb{H}_i$ denote $\mathbb{H}_i(LS^n)$.\\
\begin{tikzpicture}[descr/.style={fill=white,inner sep=1.5pt}]
        \matrix (m) [
            matrix of math nodes,
            row sep=2em,
            column sep=2.5em,
            text height=1.5ex, text depth=0.25ex
        ]
        { & H^{S^1}_{n+3} & H^{S^1}_{n+1} & \mathbb{H}_{2}\cong 0 \\        
        & H^{S^1}_{n+2} & H^{S^1}_n\cong 0 & \mathbb{H}_{1}\cong 0\\                  
             & H^{S^1}_{n+1} & H^{S^1}_{n-1} \cong \ZZ \oplus \ZZ& \mathbb{H}_{0}\cong \ZZ (1\otimes 1) \\     
        & H^{S^1}_{n}\cong 0 & H^{S^1}_{n-2}\cong 0 & \mathbb{H}_{-1}\cong \ZZ (a\otimes u) \\     
        & H^{S^1}_{n-1} & H^{S^1}_{n-3}\cong \ZZ (\gamma_{\frac{n-3}{2}}) & \mathbb{H}_{-2}\cong 0 \\     
        & H^{S^1}_{n-2}\cong 0 & H^{S^1}_{n-4} & 0\\        
        &\vdots &\vdots &\vdots\\
        & H^{S^1}_4 \cong \ZZ (\gamma_2) & H^{S^1}_2\cong \ZZ (\gamma) & \mathbb{H}_{-n+3}\cong 0 \\
        & H^{S^1}_3 & H^{S^1}_1 & \mathbb{H}_{-n+2}\cong 0 \\
          & H^{S^1}_2\cong \ZZ( \gamma) & H^{S^1}_0\cong \ZZ & \mathbb{H}_{-n+1}\cong 0 \\
          & H^{S^1}_1\cong 0 & 0 & \mathbb{H}_{-n}\cong \ZZ (a \otimes 1) \\
            & H^{S^1}_0\cong \ZZ &0.  &  \\
        };

        \path[overlay,->, font=\scriptsize,>=latex]

        (m-1-2) edge  node[above] {$\cong$}(m-1-3)
        (m-1-3) edge[orange] node[above] {$M$}(m-1-4) 
        (m-1-4) edge[out=355,in=175,cyan] node[descr,yshift=0.3ex] {$e$} (m-2-2)
        (m-2-2) edge node[above] {$\cong$} (m-2-3)
        (m-2-3) edge[orange] node[above] {$M$}(m-2-4)
        (m-2-4) edge[out=355,in=175,cyan] node[descr,yshift=0.3ex] {$e$} (m-3-2)
        (m-3-2) edge (m-3-3)
        (m-3-3) edge[orange] node[above] {$M$} (m-3-4)
        (m-3-4) edge[out=355,in=175,cyan] node[descr,yshift=0.3ex] {$e$} (m-4-2)
        (m-4-2) edge (m-4-3)
        (m-4-3) edge[orange] node[above] {$M$} (m-4-4)
         (m-4-4) edge[out=355,in=175,cyan] node[descr,yshift=0.3ex] {$e$}(m-5-2)   
        (m-5-2) edge (m-5-3)
        (m-5-3) edge[orange] node[above] {$M$}(m-5-4)
         (m-5-4) edge[out=355,in=175,cyan] node[descr,yshift=0.3ex] {$e$}(m-6-2)   
        (m-6-2) edge  node[above] {$\cong$} (m-6-3)
        (m-6-3) edge[orange] node[above] {$M$}(m-6-4)
        (m-8-2) edge node[above] {$\cong$} (m-8-3)
        (m-8-3) edge[orange] node[above] {$M$}(m-8-4)
        (m-8-4) edge[out=355,in=175,cyan] node[descr,yshift=0.3ex] {$e$}(m-9-2)   
        (m-9-2) edge node[above] {$\cong$} (m-9-3)
        (m-9-3) edge[orange] node[above] {$M$}(m-9-4)
        (m-9-4) edge[out=355,in=175,cyan] node[descr,yshift=0.3ex] {$e$}(m-10-2)   
        (m-10-2) edge node[above] {$\cong$} (m-10-3)
        (m-10-3) edge[orange] node[above] {$M$}(m-10-4)
        (m-10-4) edge[out=355,in=175,cyan] node[descr,yshift=0.3ex] {$e$}(m-11-2)   
        (m-11-2) edge node[above] {$\cong$}(m-11-3)
        (m-11-3) edge[orange] node[above] {$M$}(m-11-4)
        (m-11-4) edge[out=355,in=175,cyan] node[descr,yshift=0.3ex] {$e,\cong$}(m-12-2)   
        (m-12-2) edge (m-12-3);
\end{tikzpicture}\\
The maps $H^{S^1}_i\longrightarrow H^{S^1}_{i-2}$ are given by the cap product with the class generator $x\in H^2(\mathbb{C}P^{\infty})$. Since $H^{S^1}_2(LS^n)\cong \ZZ$, we denote the generator by $\gamma$, which is dual to $x$. We use the notation $\gamma_i=\frac{\gamma^i}{i!}$, dual to $x^i$. Since the maps given by the cap product are isomorphisms between where the loop homology is zero, we have that
\begin{align*}
H^{S^1}_{2i+1}(LS^n)=& 0, \hspace*{.5cm} 0\leq i \leq \frac{n-3}{2}\\
H^{S^1}_{2i}(LS^n)=& \ZZ (\gamma_i), \hspace*{.5cm} 1 \leq i \leq \frac{n-3}{2}.
\end{align*}
For even degrees, the generator $\gamma_i$ increases subscript as isomorphisms in the sequence are given by cap product with $x$, dual to the cup product with $x$.\\
\indent To determine $H^{S^1}_{n-1}(LS^n)$, note that we have a short exact sequence,
\begin{equation}
0\longrightarrow \ZZ (a \otimes u) \longrightarrow H^{S^1}_{n-1} \longrightarrow  \ZZ (\gamma_{\frac{n-3}{2}}) \longrightarrow 0\nonumber
\end{equation}
that splits since the last term is free. Thus $H^{S^1}_{n-1}(LS^n)\cong \ZZ e(a \otimes u) \oplus (\gamma_{\frac{n-1}{2}})$. We use the notation of $e(-)$ to denote that the generator comes from the erasing map. Using the BV-operator to determine the marking map $M:H^{S^1}_{n-1} \rightarrow \mathbb{H}_0$, we have that $M(a \otimes u)=1\otimes 1$, so the erasing map $e:\mathbb{H}_0\rightarrow H^{S^1}_n$ is zero, thus $H_n^{S^1}(LS^n)\cong 0$. The marking map sends generators coming from $H_*(\mathbb{C}P^{\infty})$ to zero. 
\lem \label{lem:marking} $M(\gamma_{\frac{n-1}{2}})=0$.
\proof 
In the circle bundle~(\ref{eqn:oddcirclebundle})
the marking map is an umkehr map coming from the projection map. Notice that $\mathbb{C}P^{\infty}=BS^1=\{pt\}\times ES^1 \subset LS^n \times_{S^1} ES^1$. Since $\pi^{-1}(\{pt\}\times_{S^1} ES^1)=\{pt\} \times ES^1$, which is contractible, then $M$ maps generators from $\mathbb{C}P^{\infty}$ into a contractible space, thus $M(\gamma_{i})=0$ for any $i$, where $\gamma_i$ denotes a generator coming from the homology of $\mathbb{C}P^{\infty}$.
\qed \\
\normalfont
\indent With the knowledge that $M(\gamma_{\frac{n-1}{2}})=0$, the cap product map $H^{S^1}_{n+1}\rightarrow H^{S^1}_{n-1}$ is injective with image isomorphic to $\ZZ (\gamma_{\frac{n-1}{2}})$ so we have $H^{S^1}_{n+1}\cong \ZZ (\gamma_{\frac{n+1}{2}})$.\\
\indent In the next piece of the Gysin sequence where loop homology is non-zero, we have\\
\begin{center}
\begin{tikzpicture}[descr/.style={fill=white,inner sep=1.5pt}]
        \matrix (m) [
            matrix of math nodes,
            row sep=2em,
            column sep=2.5em,
            text height=1.5ex, text depth=0.25ex
        ]
        {         
         H^{S^1}_{2n+1} & H^{S^1}_{2n-1} & \mathbb{H}_{n}\cong 0\\                  
              H^{S^1}_{2n} & H^{S^1}_{2n-2} \cong \ZZ \oplus \ZZ& \mathbb{H}_{n-1}\cong \ZZ (1\otimes u) \\     
         H^{S^1}_{2n-1} & H^{S^1}_{2n-3}\cong 0 & \mathbb{H}_{n-2}\cong \ZZ (a\otimes u^2) \\     
         H^{S^1}_{2n-2} & H^{S^1}_{2n-4}\cong \ZZ (\gamma_{n-2}) & \mathbb{H}_{n-3}\cong 0 \\     
         H^{S^1}_{2n-3}\cong 0 & H^{S^1}_{2n-5}\cong 0 & 0\\        
        \vdots &\vdots &\vdots\\
        H^{S^1}_{n+3} \cong \ZZ (\gamma_{n-2}) & H^{S^1}_{n+1}\cong \ZZ (\gamma_{\frac{n+1}{2}}) & \mathbb{H}_{2}\cong 0 \\
         H^{S^1}_{n+2} & H^{S^1}_n & \mathbb{H}_{1}\cong 0. \\
        };

        \path[overlay,->, font=\scriptsize,>=latex]

        (m-1-1) edge  node[above] {$\cong$}(m-1-2)
        (m-1-2) edge[orange] node[above] {$M$}(m-1-3) 
        (m-1-3) edge[out=355,in=175,cyan] node[descr,yshift=0.3ex] {$e$} (m-2-1)
        (m-2-1) edge  (m-2-2)
        (m-2-2) edge[orange] node[above] {$M$}(m-2-3)
        (m-2-3) edge[out=355,in=175,cyan] node[descr,yshift=0.3ex] {$e$} (m-3-1)
        (m-3-1) edge (m-3-2)
        (m-3-2) edge[orange] node[above] {$M$} (m-3-3)
        (m-3-3) edge[out=355,in=175,cyan] node[descr,yshift=0.3ex] {$e$} (m-4-1)
        (m-4-1) edge (m-4-2)
        (m-4-2) edge[orange] node[above] {$M$} (m-4-3)
         (m-4-3) edge[out=355,in=175,cyan] node[descr,yshift=0.3ex] {$e$}(m-5-1)   
        (m-5-1) edge (m-5-2)
        (m-5-2) edge[orange] node[above] {$M$}(m-5-3)

                (m-7-1) edge node[above] {$\cong$} (m-7-2)
        (m-7-2) edge[orange] node[above] {$M$}(m-7-3)
         (m-7-3) edge[out=355,in=175,cyan] node[descr,yshift=0.3ex] {$e$}(m-8-1)   
        
        (m-8-1) edge node[above] {$\cong$} (m-8-2)
        (m-8-2) edge[orange] node[above] {$M$}(m-8-3)

    ;
\end{tikzpicture}\\
\end{center}
In the third and fourth row above, we have a short exact sequence with $H^{S^1}_{2n-2}$ in the center, which splits, so $H^{S^1}_{2n-2}\cong \ZZ (e(a\otimes u^2)) \oplus \ZZ (\gamma_{n-1})$. Mapping $H^{S^1}_{2n-2}$ through $M$, we have $M(e(a\otimes u^2))=2(1\otimes u)$, given by the BV-operator.
\indent Thus $H^{S^1}_{2n-1}\cong \ZZ/2\ZZ (e(1 \otimes u))$. Since the cap product map $H^{S^1}_{2n}\rightarrow H^{S^1}_{2n-2}$ is injective with image $\ZZ (\gamma_{n-1})$, $H^{S^1}_{2n} \cong \ZZ (\gamma_{n})$. Summarizing, we have
\begin{align*}
H^{S^1}_{2i+1}(LS^n)\cong & 0, \hspace*{.5cm} & \frac{n-1}{2}\leq i \leq n-2 \\
H^{S^1}_{2i}(LS^n)\cong & \ZZ (\gamma_i), \hspace*{.5cm} & \frac{n+1}{2} \leq i \leq n-2  \\
H^{S^1}_{n-2}(LS^n)\cong & \ZZ e(a\otimes u^2) \oplus \ZZ (\gamma_{n-1}) \\
H^{S^1}_{n-1}(LS^n)\cong & \ZZ/2\ZZ (e(1 \otimes u)) \\
H^{S^1}_{2n}(LS^n)\cong & \ZZ (\gamma_{n}) \\
H^{S^1}_{2i+1} \cong & \ZZ_2 (1 \otimes u)\gamma_{i-n-1} \hspace*{.5cm}& n-1 \leq i \leq \frac{3n-5}{2} \\
H^{S^1}_{2i} \cong & \ZZ \gamma_i \hspace*{.5cm} & \frac{n}{2} \leq i \leq \frac{3n-5}{2} 
\end{align*}
Now assume the following holds for all $k\in \NN$:

\begin{equation} \label{eqn:oddinduction}
H^{S^1}_{2i} \cong \left\{
\begin{array}{rl}
\ZZ \gamma_i & \text{if } (n-1)\nmid 2i,\\
\ZZ \gamma_i \oplus \ZZ(e(a\otimes u^i)) & \text{if } (n-1)\vert 2i.
\end{array} \right.
\end{equation}
for $1\leq i \leq \frac{(k+1)(n-1)-2}{2}$ and 
\begin{equation}
H_{2i+1}^{S^1}\cong t_k \nonumber
\end{equation}
for $\frac{k(n-1)}{2}\leq i \leq \frac{(k+1)(n-1)-2}{2}$, where $t_k$ is a torsion group of order $k!$. There are extension issues, so we cannot say which torsion group $H_{2i+1}^{S^1}$ should be.)\\
\indent Consider the $k+1$-th non-zero piece of the Gysin sequence:\\
\vspace{.2cm}
\begin{center}
\begin{tikzpicture}
 [descr/.style={fill=white,inner sep=1.5pt}]
        \matrix (m) [
            matrix of math nodes,
            row sep=2em,
            column sep=.5em,
            text height=1.5ex, text depth=0.25ex
        ]
        {         
        & H^{S^1}_{(k+1)(n-1)+4} & H^{S^1}_{(k+1)(n-1)+2} & 0\\                  
             & H^{S^1}_{(k+1)(n-1)+3} & H^{S^1}_{(k+1)(n-1)+1}& 0\\     
        & H^{S^1}_{(k+1)(n-1)+2} & H^{S^1}_{(k+1)(n-1)} &  \ZZ (1\otimes u^k) \\     
        & H^{S^1}_{(k+1)(n-1)+1} & H^{S^1}_{(k+1)(n-1)-1}\cong \bigoplus_{j=2}^{k}\ZZ_j &\ZZ (a \otimes u^{k+1}) \\     
        & H^{S^1}_{(k+1)(n-1)} & H^{S^1}_{(k+1)(n-1)-2}\cong \ZZ(\gamma_{\frac{(k+1)(n-1)-2}{2}}) & 0.\\    
        };

        \path[overlay,->, font=\scriptsize,>=latex]

        (m-1-2) edge  node[above] {$\cong$}(m-1-3)
        (m-1-3) edge[orange] node[above] {$M$}(m-1-4) 
        (m-1-4) edge[out=355,in=175,cyan] node[descr,yshift=0.3ex] {$e$} (m-2-2)
        (m-2-2) edge  (m-2-3)
        (m-2-3) edge[orange] node[above] {$M$}(m-2-4)
        (m-2-4) edge[out=355,in=175,cyan] node[descr,yshift=0.3ex] {$e$} (m-3-2)
        (m-3-2) edge (m-3-3)
        (m-3-3) edge[orange] node[above] {$M$} (m-3-4)
        (m-3-4) edge[out=355,in=175,cyan] node[descr,yshift=0.3ex] {$e$} (m-4-2)
        (m-4-2) edge (m-4-3)
        (m-4-3) edge[orange] node[above] {$M$} (m-4-4)
         (m-4-4) edge[out=355,in=175,cyan] node[descr,yshift=0.3ex] {$e$}(m-5-2)   
        (m-5-2) edge (m-5-3)
        (m-5-3) edge[orange] node[above] {$M$}(m-5-4)

    ;
\end{tikzpicture}\\
\end{center}
\normalfont
Thus, we can extract a short exact sequence from the last two lines of the Gysin sequence above, giving us $H^{S^1}_{(k+1)(n-1)}\cong \ZZ(e(a\otimes u^{k+1}))\oplus \ZZ (\gamma_{\frac{(k+1)(n-1)}{2}})$.
\normalfont
\indent It can be seen that $H^{S^1}_{(k+1)(n-1)+2}\cong \ZZ (\gamma_{\frac{(k+1)(n-1)+2}{2}})$ and that\\ $torsion(H^{S^1}_{(k+1)(n-1)+1})\cong t_{k+1}$. $H^{S^1}_{(k+1)(n-1)+1}$ is all torsion since it is sandwiched between a short exact sequence of torsion groups, so $H^{S^1}_{(k+1)(n-1)+1}\cong  t_{k+1}$. Since loop homology $\mathbb{H}_{i}(LS^n)$ is zero for $(k+1)(n-1)+2-n \leq i \leq (k+2)(n-1)-1-n$), we obtain the analogous statements of~(\ref{eqn:oddinduction}) for $k+1$.
\qed \\
\indent The string bracket is a degree $2-n$ map, where the only possible non-zero bracket is of the generators $e(a\otimes u^i)$, since the marking map $M$ sends all other generators to zero.
\begin{proof}[Proof of Theorem \ref{thm:stringbracket} (2)]
\begin{align*}
[e(a\otimes u^i), e(a\otimes u^j)]=&e(M(e(a\otimes u^i))\bullet M(e(a\otimes u^j)) \nonumber \\
=&e(i(1\otimes u^{i-1})\bullet j(1\otimes u^{j-1})) \nonumber \\
=&e(ij(1\otimes u^{i+j-2})). \nonumber 
\end{align*}
where $e(1\otimes u^{i+j-2})$ is a generator of $\ZZ_{i+j-1}$, so the bracket is only zero when $i+j-1$ divides $ij$.
\end{proof}
\subsection{String Homology for Even Spheres}
For $n$ even,
\begin{align*}
\mathbb{H}_*(LS^{n},\ZZ) \cong \frac{\Lambda(b) \otimes \ZZ[a,v]}{(a^2,ab,2av)}
\end{align*}
where $|a|=-n$, $|b|=-1$ and $|v|=2n-2$ \cite{CJY}. By Menichi, we have that the BV-operator acts as follows \cite{Me}:
\begin{align*}
\Delta(v^k) &=0\\
\Delta(av^k) &=0\\
\Delta(bv^k) &= (2k+1)v^k.
\end{align*}
To keep track of things, $|av^k|=k(2n-2)-n$, $|bv^k|=k(2n-2)-1$, $|v^k|=k(2n-2)$. Let us consider the bottom of the Gysin sequence:\\

\begin{tikzpicture}[descr/.style={fill=white,inner sep=1.5pt}]
        \matrix (m) [
            matrix of math nodes,
            row sep=2em,
            column sep=2.5em,
            text height=1.5ex, text depth=0.25ex
        ]
        {       & H^{S^1}_{3n-1} & H^{S^1}_{3n-2} & \mathbb{H}_{2n-2}\cong \ZZ(v^2) \\                  
             & H^{S^1}_{3n-2} & H^{S^1}_{3n-4}& \mathbb{H}_{2n-3}\cong \ZZ (bv) \\   
             &\vdots &\vdots &\vdots\\  
        & H^{S^1}_{2n-1} & H^{S^1}_{2n-3}\cong 0 & \mathbb{H}_{n-2}\cong \ZZ_2 (av) \\     
        &\vdots &\vdots &\vdots\\
        & H^{S^1}_{n+1} & H^{S^1}_{n-1}\cong \ZZ & \mathbb{H}_{0}\cong \ZZ (v) \\     
        & H^{S^1}_{n} & H^{S^1}_{n-2} & \mathbb{H}_{-1}\cong \ZZ (b)\\        
        &\vdots &\vdots &\vdots\\
          & H^{S^1}_2 & H^{S^1}_0 & \mathbb{H}_{-n+1}\cong 0 \\
          & H^{S^1}_1 & 0 & \mathbb{H}_{-n}\cong \ZZ (a) \\
            & H^{S^1}_0 &0  &  \\
        };

        \path[overlay,->, font=\scriptsize,>=latex]

        (m-1-2) edge  node[above] {$\cong$}(m-1-3)
        (m-1-3) edge[orange] node[above] {$M$}(m-1-4) 
        (m-1-4) edge[out=355,in=175,cyan] node[descr,yshift=0.3ex] {$e$} (m-2-2)
        (m-2-2) edge node[above] {$\cong$} (m-2-3)
        (m-2-3) edge[orange] node[above] {$M$}(m-2-4)
        
        (m-4-2) edge (m-4-3)
        (m-4-3) edge[orange] node[above] {$M$} (m-4-4)
         
        (m-6-2) edge (m-6-3)
        (m-6-3) edge[orange] node[above] {$M$}(m-6-4)
         (m-6-4) edge[out=355,in=175,cyan] node[descr,yshift=0.3ex] {$e$}(m-7-2)   
        (m-7-2) edge  node[above] {$\cong$} (m-7-3)
        (m-7-3) edge[orange] node[above] {$M$}(m-7-4)
 (m-9-2)   
        (m-9-2) edge node[above] {$\cong$} (m-9-3)
        (m-9-3) edge[orange] node[above] {$M$}(m-9-4)
        (m-9-4) edge[out=355,in=175,cyan] node[descr,yshift=0.3ex] {$e$}(m-10-2)   
        (m-10-2) edge node[above] {$\cong$}(m-10-3)
        (m-10-3) edge[orange] node[above] {$M$}(m-10-4)
        (m-10-4) edge[out=355,in=175,cyan] node[descr,yshift=0.3ex] {$e,\cong$}(m-11-2)   
        (m-11-2) edge (m-11-3);
\end{tikzpicture}\\
from this sequence and knowledge of the BV-operator, we get
\pagebreak
\begin{align*}
H_0^{S^1}(LS^n)&\cong \ZZ a\\
H_1^{S^1}(LS^n)&\cong 0\\
H_2^{S^1}(LS^n)&\cong \ZZ \gamma \\
H_3^{S^1}(LS^n)&\cong 0 \\
H_4^{S^1}(LS^n)&\cong \ZZ \gamma_2\\
\vdots &\\
H_{n-1}^{S^1}(LS^n)&\cong \ZZ e(b)\\
H_{n}^{S^1}(LS^n)&\cong \ZZ e(v) \oplus \ZZ \gamma_{\frac{n}{2}} \\
\vdots &\\
H_{2n-2}^{S^1}(LS^n)&\cong \ZZ_2 e(av) \oplus \ZZ \oplus \ZZ \gamma_{\frac{2n-2}{2}} \\
H_{2n-1}^{S^1}(LS^n)&\cong 0 \\
H_{2n}^{S^1}(LS^n)&\cong \ZZ_2  \oplus \ZZ \oplus \ZZ \\
\vdots & \\
H_{3n-3}^{S^1}(LS^n)&\cong \ZZ e(bv)\\
H_{3n-2}^{S^1}(LS^n)&\cong \ZZ_2  \oplus \ZZ_3 e(v^3) \oplus \ZZ \oplus \ZZ \\
\vdots & \\
H_{4n-4}^{S^1}(LS^n)&\cong \ZZ_2 \oplus \ZZ_2 \oplus \ZZ_3 \oplus \ZZ \oplus \ZZ \\
H_{4n-3}^{S^1}(LS^n)&\cong 0\\
\vdots \\
H_{5n-5}^{S^1}(LS^n)&\cong \ZZ e(bv^2)\\
H_{5n-4}^{S^1}(LS^n)&\cong \ZZ_2 \oplus \ZZ_2 \oplus \ZZ_2 \oplus \ZZ_3 \oplus \ZZ \oplus \ZZ \\
\vdots &
\end{align*}
where all of the odd degree homology are isomorphic, and all even degree homology are isomorphic, or $H_i^{S^1}\cong H_{i-2}^{S^1}$ in the gaps denoted by the vertical dots. The $k$-th piece of the sequence is as follows:

\begin{tikzpicture}[descr/.style={fill=white,inner sep=1.5pt}]
        \matrix (m) [
            matrix of math nodes,
            row sep=2em,
            column sep=2.5em,
            text height=1.5ex, text depth=0.25ex
        ]
        {       & H^{S^1}_{k(2n-2)+n+1} & H^{S^1}_{k(2n-2)+n-1} & \mathbb{H}_{k(2n-2)}\cong \ZZ(v^k) \\                  
             & H^{S^1}_{k(2n-2)+n} & H^{S^1}_{k(2n-2)-2+n}& \mathbb{H}_{k(2n-2)-1}\cong \ZZ (bv^k) \\   
             &\vdots &\vdots &\vdots\\  
        & H^{S^1}_{k(2n-2)+1} & H^{S^1}_{k(2n-2)-1}\cong 0 & \mathbb{H}_{k(2n-2)-n}\cong \ZZ_2 (av^k) \\     
        };

        \path[overlay,->, font=\scriptsize,>=latex]

        (m-1-2) edge  node[above] {$\cong$}(m-1-3)
        (m-1-3) edge[orange] node[above] {$M$}(m-1-4) 
        (m-1-4) edge[out=355,in=175,cyan] node[descr,yshift=0.3ex] {$e$} (m-2-2)
        (m-2-2) edge node[above] {$\cong$} (m-2-3)
        (m-2-3) edge[orange] node[above] {$M$}(m-2-4)
        
        (m-4-2) edge (m-4-3)
        (m-4-3) edge[orange] node[above] {$M$} (m-4-4)
         
       ;
\end{tikzpicture}\\
inductively, we have that 
\begin{align*}
H_{k(2n-2)-2}^{S^1} \cong \ZZ_2^{k-1}\oplus C_k \oplus \ZZ \oplus \ZZ
\end{align*} 
where $C_k$ is a torsion group of order  $\prod_{i=1}^{k-1}(2i+1)$. The bottom of the above Gysin sequence gives the short exact sequence
\begin{align*}
0 \rightarrow \ZZ_2 av^k \rightarrow H^{S^1}_{k(2n-2)} \rightarrow H_{k(2n-2)-2}^{S^1} \rightarrow 0
\end{align*} 
which gives
\begin{align*}
H^{S^1}_{k(2n-2)}\cong \ZZ_2 e(av^k) \oplus \ZZ_2^{k-1}\oplus \left(\text{torsion group of order } \sum_{i=1}^{k-1}(2i+1)\right) \oplus \ZZ \oplus \ZZ.
\end{align*}
Note that the even torsion can be resolved using the results by Westerland in \cite{We}. From the top of the above Gysin sequence, we get the following.
\begin{align*}
H_{k(2n-2)+1}^{S^1}(LS^n)& \cong 0\\
\vdots & \\
H_{k(2n-2)+n-1}^{S^1}(LS^n) & \cong \ZZ e(bv^k)\\
H_{k(2n-2)+n}^{S^1}(LS^n) & \cong \ZZ_2^{k} \oplus C_k \oplus \ZZ \oplus \ZZ . 
\qed
\end{align*}
The string bracket is always zero except on the generators $e(bv^j)$ since $e(v^{k+l})$ has order $2(k+l)+1$ so it is not always zero.
\begin{proof}[Proof of Theorem \ref{thm:stringbracket} (3)]
\begin{align*}
[e(bv^k),e(bv^l)]&=(-1)^{k(2n-2)-1-n}e(M(e(bv^k))\bullet M(e(bv^l))\\
&=-(4kl+2k+2l+1)e(v^{k+l}).
\end{align*}
\end{proof}

\section{Examples Computation for String Homology and String Bracket of $S^1$.}
\indent We compute the string homology of $S^1$ using the Gysin sequence for the circle bundle
\begin{eqnarray}
 \nonumber S^1 \rightarrow LS^1 \times ES^1 \rightarrow LS^1 \times_{S^1} ES^1
\end{eqnarray}\\
Basu computed this in his thesis, but here we use different techniques. \\
\begin{proof}[Proof of Theorem \ref{thm:stringhomology1} (1)]
\indent Recall that the non-equivariant homology of $LS^1$ is given as follows \cite{CJY}, \cite{He}, 
\begin{equation}
\nonumber \mathbb{H}_*(LS^1)=\Lambda_{\ZZ}[a]\otimes \ZZ[x,x^{-1}], \; |a |=-1, |x|=0,
\end{equation}
where $\mathbb{H}_*(LS^1)=H_{*+1}(LS^1)$ and $a$ corresponds to the dual of  $[S^1]$ under the geometric grading \cite{Se}, \cite{CJY}.\\
\indent The BV-operator ($\Delta=M\circ e$) acts on generators of $\mathbb{H}_*(LS^1)$ as follows, \cite{Me}: 
\begin{eqnarray}
&\Delta (a \otimes x^i)& =i(1\otimes x^{i})\nonumber \\
&\Delta (1\otimes x^i) & =0.\nonumber 
\end{eqnarray}
\indent Consider the Gysin sequence for the above circle bundle:\\
\begin{tikzpicture}[descr/.style={fill=white,inner sep=1.5pt}]
        \matrix (m) [
            matrix of math nodes,
            row sep=2em,
            column sep=2.5em,
            text height=1.5ex, text depth=0.25ex
        ]
        { &  & &0 \\
          & H^{S^1}_2(LS^1) & H^{S^1}_0(LS^1) & \mathbb{H}_{0}(LS^1)\cong \bigoplus_{n\in \ZZ} \ZZ (1\otimes x^n) \\
            & H^{S^1}_1(LS^1) & H^{S^1}_{-1}(LS^1)\cong 0 & \mathbb{H}_{-1}(LS^1)\cong \bigoplus_{n\in \ZZ} \ZZ (a\otimes x^n) \\
            & H^{S^1}_0(LS^1) & 0. & \\
           \\
        };
        \path[overlay,->, font=\scriptsize,>=latex]
    
        (m-1-4) edge[out=355,in=175,cyan] node[descr,yshift=0.3ex] {$e$} (m-2-2)
        (m-2-2) edge  node[above] {$c$} (m-2-3)
        (m-2-3) edge[orange] node[above] {$M$}(m-2-4)
        (m-2-4) edge[out=355,in=175,cyan] node[descr,yshift=0.3ex] {$e$} (m-3-2)
        (m-3-2) edge (m-3-3)
        (m-3-3) edge[orange] node[above] {$M$} (m-3-4)
        (m-3-4) edge[out=355,in=175,cyan] node[descr,yshift=0.3ex] {$e$} (m-4-2)
        (m-4-2) edge (m-4-3)
         ;
\end{tikzpicture}\\
\indent The end of the Gysin sequence gives us that ${H^{S^1}_0(LS^1)\cong \bigoplus\limits_{n\in \ZZ} \ZZ (e(a\otimes x^n))}$. Using the information from the BV-operator, ${M\circ e(a \otimes x^n)=\Delta (a \otimes x^n)}=n(1 \otimes x^n)$. Since $e$ is surjective and ${ker(e)=im(M)= \bigoplus\limits_{n\in \ZZ} n\ZZ(1 \otimes x^n)}$, we have that ${H^{S^1}_1(LS^1) \cong  \mathbb{H}_{0}(LS^1) / ker(e) \cong \bigoplus\limits_{n\in \ZZ} \ZZ /n\ZZ (1 \otimes x^n) \oplus \ZZ (1\otimes 1)} $. At the beginning of the Gysin sequence, we have ${im(c)=ker(M)=\ZZ (a\otimes 1)}$, and since $c$ is injective, ${H^{S^1}_2(LS^1) \cong \ZZ (a \otimes 1)}$. Summarizing, we get the following.
\begin{eqnarray}
& H^{S^1}_0(LS^1) &\cong \bigoplus_{n\in \ZZ} \ZZ (e(a\otimes x^n))\nonumber \\
& H^{S^1}_1(LS^1) &\cong  H^{S^1}_{2i+1}(LS^1) \cong \bigoplus_{n\in \ZZ -\{0\}} \ZZ /n\ZZ (1 \otimes x^n)\oplus \ZZ (1 \otimes 1), \hspace*{.5cm} i\geq 0\nonumber \\
& H^{S^1}_2(LS^1) &\cong H^{S^1}_{2i}(LS^1) \cong \ZZ (a \otimes 1), \hspace*{.5cm} i\geq 1.\nonumber 
\end{eqnarray}\\
\end{proof}
\normalfont
\begin{proof}[Proof of Theorem \ref{thm:stringbracket} (1)]
The string bracket, 
\begin{align*}
H_i^{S^1}(LS^1) \otimes H_j^{S^1}(LS^1) \rightarrow H_{i+j+1}^{S^1}(LS^1)
\end{align*}  is a degree $+1$ map, and it is only nontrivial on generators of degree zero since the marking map is trivial for generators of degree greater than zero. For ${e(a \otimes x^n)}$, ${e(a \otimes x^m)}$ in ${H^{S^1}_0(LS^1)}$,
\begin{align*}
[e(a \otimes x^n),e( a \otimes x^m)] & =(-1)^{-1}e(M(e(a \otimes x^n))\bullet M(e(a \otimes x^m))) \\
&=-e(n(1 \otimes x^n) \bullet m(1 \otimes x^m))\\
& =-nm(e(1 \otimes x^{n+m}))\\
&=-nm(1 \otimes x^{n+m}).
\end{align*}
So ${[a \otimes x^n, a \otimes x^m]=0}$ if $n+m \neq 0$ and $n+m$ divides $nm$. If $n+m=0$ then  ${[a \otimes x^n, a \otimes x^m]=nm(1\otimes 1)}$.
We can conclude that the bracket is only nontrivial for the torsion elements. 
\end{proof}

\bibliographystyle{plain}
\bibliography{bibfile}

\end{document}